\documentclass[aps, 11pt] {revtex4}
\usepackage[utf8]{inputenc}
\usepackage[english]{babel}
\usepackage[titletoc,toc,title]{appendix}
\usepackage{amsmath}
\usepackage{amsfonts}
\usepackage{amssymb}
\usepackage[colorlinks=true,linkcolor=blue,citecolor=red,urlcolor=blue]{hyperref}
\usepackage{graphicx}
\usepackage{enumerate}
\usepackage{color}
\usepackage{mathtools}
\usepackage{hyperref}

\usepackage{subfigure}

\newcommand{\ms}{\medskip}

\begin{document}

	\title{Isochronous waveforms of Li\'enard equations via commutative factorization}
	
	\author{G. Gonz\'alez}
\email{ gabriel.gonzalez@academicos.udg.mx, ORCID: 0000-0001-8217-321x} 

	\affiliation{
Departamento de Ciencias B\'asicas y Aplicadas, CUTonal\'a, Universidad de Guadalajara, Avenida Nuevo Perif\'erico 555, Ejido San Jos\'e Tateposco, 45425 Tonal\'a, Mexico}

\author{O. Cornejo-P\'erez}
\email{octavio.cornejo@uaq.mx, ORCID: 0000-0002-1790-7640}
\affiliation{Facultad de Ingenier\'{\i}a, Universidad Aut\'onoma de Quer\'etaro,\\
Centro Universitario Cerro de las Campanas, 76010 Santiago de Quer\'etaro, Qro., Mexico}

\author{J. de la Cruz}
\email{josue.dl@slp.tecnm.mx, 
ORCID: 0000-0001-5943-5752}
\affiliation{TECNM-Instituto Tecnol\'ogico de San Luis Potos\'{\i}, \\
Avenida Tecnol\'ogico S/N Col. Unidad Ponciano Arriaga, 78437 Soledad de Graciano S\'anchez, S.L.P., Mexico}

\author{H.~C. Rosu}
\email{hcr@ipicyt.edu.mx, ORCID: 0000-0001-5909-1945, Corresponding author.}
\affiliation{IPICyT, Instituto Potosino de Investigaci\'on Cient\'{\i}fica y Tecnol\'ogica,\\
Camino a la presa San Jos\'e 2055, Col. Lomas 4a Secci\'on, 78216 San Luis Potos\'{\i}, S.L.P., Mexico}


\begin{abstract}
\textbf {Abstract.} 
{\scriptsize  Isochronous waveform solutions of homogeneous Li\'enard equations
are obtained by a modification of the nonlinear factorization method of Rosu and Cornejo-P\'erez. 
 The scheme is based on the assumption that the intermediate function $\Phi$ that can be introduced in this factorization method depends on both 
 the dependent and independent variables of the nonlinear equation.
 The method is applied to three cases, a noted cubic anharmonic oscillator, a Li\'enard-reduced form of the Sharma-Tasso-Olver evolution equation, and the cubic-quintic Wilson's Li\'enard equation. All these cases are written in a commutative factored form that allows to obtain the general solutions as solutions of a certain type of Bernoulli differential equation. A theorem is also given asserting the general form of the Li\'enard equation, i.e., for given polynomial degree $n$ of its coefficients, which can be solved by this method. The conditions under which these equations can be also approached by non-local transformations are established.}

\medskip

\end{abstract}

\centerline{Phys. Lett. A 564, 131087 (2025)}
\centerline{10.1016/j.physleta.2025.131087}

\vspace{2pc}

\maketitle

\newpage

\section*{1. Introduction}\label{sec1}

The nonlinear ordinary differential equations widely known as of the Li\'enard type \cite{lrbook,tiwari2015,hl16,demina2023} have been introduced in 1928 by Li\'enard in two papers written in French, in which he studied the generalized van der Pol equation
\begin{equation}\label{eq00}
\frac{d^2x}{dt^2} +\omega f(x)\frac{dx}{dt} + \omega^2x = 0~.
\end{equation}
In 1934, in a brief overview of many of his papers \cite{lienard1934}, Li\'enard also mentions the equation 
\begin{equation}\label{eq00}
\frac{d^2x}{dt^2} +\omega f(x)\frac{dx}{dt} + \phi(x) = 0~,
\end{equation}
but he never published the results he seemed to have.

 On the other hand, in 1942, Levinson and Smith \cite{ls42} were the first to study such equations in their standard form
\begin{equation}\label{eq1}
\ddot{x} +G(x)\dot{x} + F(x) = 0~,
\end{equation}
and also the more general equations with $G(x,\dot{x})$, where the dot denotes the derivative $D=d/dt$ and $G(x)$ and $F(x)$ are arbitrary smooth real functions. Nowadays, the Li\'enard equations and their generalizations serve as an important
class of nonlinear ordinary differential equations with widespread applications to a huge number of physical, chemical, biological, and engineering systems, see \cite{getal} and references therein. The linear differential equations correspond to the particular Li\'enard case $F(x)\sim x$ and $G(x)=$ constant.
In the literature, the integrability of the Li\'enard equations was naturally an important issue that has been dealt with by a variety of approaches, 
such as the Prelle-Singer method \cite{csl2005-1}, Lie point symmetries \cite{p2009-1,p2009-2}, and reductions to integrable Abel equations  \cite{mr2013,hlm2014} and Painl\'eve type equations \cite{ks2016,s2020}, among others.

\medskip

The goal of this paper is to present a novel approach to the nonlinear factorization method of these equations introduced in \cite{rcp1,rcp2} 
by which one can obtain in a straightforward way isochronous (periodic) solutions
of many Li\'enard equations with polynomial nonlinearities.

\medskip

The novelty is based on assuming that an intermediate function $\Phi$ that can be introduced
in the commutative factorization setting \cite{getal} depends on both the dependent and independent variables, $x$ and $t$, respectively.
This leads to a quasi-linear differential equation for $\Phi$, similar to the so-called field method introduced by
Vujanovic \cite{vuj1,vuj2}, see also \cite{kov1,kov2}.
The implementation of this approach in the commutative factoring framework is presented in Section 2 of the paper.
In Section 3, the way the method works is demonstrated for some specific examples. In Section 4, a theorem is provided that defines the polynomial 
Li\'enard equations of given order $n$ that can be approached by this method.
Finally, the conclusions are summarized in Section 5.

\section*{2. Commutative Nonlinear Factorization} 

Equation~(\ref{eq1}) can be factored in the form
\begin{equation}\label{eq2}
[ D-\phi_2(x) ] [ D-\phi_1(x) ]x=0~,
\end{equation}
under the conditions \cite{rcp1,rcp2}
\begin{eqnarray}
&\qquad \quad \phi_2+\frac{d(\phi_1x)}{dx}=-G(x)~, \label{eq3}\\
&\phi_1 \phi_2 x =F(x)~.\label{eq3bis}
\end{eqnarray}
From (\ref{eq3bis}), one can see that if $F(x)$ is a homogeneous polynomial then one can choose the factorization functions $\phi_1, \phi_2$ as complementary
pairs in the set of divisors of the nonhomogeneous polynomial $F(x)/x$, a hint that reveals the efficiency of this factorization method in the case of polynomial nonlinearities.
The commutative factorization setting can be achieved by demanding the following conditions
\begin{align}\label{eqcc}
  & \phi_1+\phi_2+x\frac{d\phi_1}{dx}= -G(x)~, \\
  & \phi_1+\phi_2+x\frac{d\phi_2}{dx}=-G(x)~, \\
 & \phi_1\phi_2=\frac{F(x)}{x}~,
\end{align}
which implies that the two factorization functions differ only by a constant, i.e. $\phi_2=\phi_1+C$.
\medskip

Rosu and Cornejo-P\'erez considered initially only the solutions provided by $[D-\phi_1(x)]x=0$, but also the extension to $[D-\phi_1(x)]x=\Phi(x)$ has been discussed in the literature \cite{wl2008,rmex17}. Borrowing from the field method, we assume
now $[D-\phi_1(x)]x=\Phi(x,t)$, which yields the following coupled ODE's for the factorized Eq.~(\ref{eq2}),
\begin{eqnarray}
&\qquad  \dot{x}-\phi_1(x)x=\Phi(x,t)~, \label{eq5}\\
&\dot{\Phi}-\phi_2(x)\Phi =0~. \label{eq4}
\end{eqnarray}
According to the chain rule, the system of equations~(\ref{eq5})-(\ref{eq4}) can be rewritten as the quasi-linear first order partial differential equation
\begin{equation}
\frac{\partial\Phi}{\partial x} \dot{x} + \frac{\partial\Phi}{\partial t} \equiv \frac{\partial\Phi}{\partial x} \left( \Phi + x\phi_1  \right)+ \frac{\partial\Phi}{\partial t} = \phi_2 \Phi~, \label{eq6}
\end{equation}
which can be solved by proposing the {\em ansatz} $\Phi(x,t)=x\zeta(t)$ \cite{vuj1,kovacic2003}. Then, for the function $\zeta$, one obtains the first-order differential equation
\begin{equation}\label{eq7}
\frac{d\zeta(t)}{dt} + \zeta^2(t) = (\phi_2-\phi_1) \zeta(t)~.
\end{equation}
If the commutative factorization condition is considered, then the subtraction of factorization functions $\phi_2 -\phi_1 = -2c$ holds
for $c\equiv$ const., and (\ref{eq7}) reduces to
\begin{equation} \label{eq8}
\frac{d\zeta(t)}{dt} + 2c \zeta(t)=- \zeta^2(t)~,
\end{equation}
which is a Bernoulli differential equation of nonlinear order two. Its solution is given as follows
\begin{eqnarray}
\zeta_0(t)&=&\frac{1}{t-t_0} \quad \textrm{for} \,\, c=0~, \label{eq9a}\\
\zeta_h(t) &=& -c \left[1- \tanh(ct+\delta)\right] \quad \textrm{for} \,\, c \in \mathbb{R} \setminus \{0\}~, \label{eq9b}\\
\zeta_t(t) &=& -\tilde{c}  \left[i+\tan( \tilde{c} t+\delta)\right] \quad \textrm{for} \,\, c =i\tilde{c} \in \mathbb{I} \setminus \{i0\}~, \label{eq9c}
\end{eqnarray}
where $\delta$ is an integration constant.
This last result allows to rewrite the first-order ODE (\ref{eq5}) in the form
\begin{equation}\label{eq10}
\dot{x}-\zeta(t)x=\phi_1(x)x~,
\end{equation}
which is a Bernoulli differential equation with the nonlinearity of one order higher than the order of the $\phi$'s.
Its solution provides the general solution of Eq.~(\ref{eq1}) factored in the (commutative) form given in Eq.~(\ref{eq2}).

It is worth mentioning that for the complex case $\phi_2 -\phi_1 = -2i\tilde{c}$, the following symmetric factorization of Eq.~(\ref{eq1})
is necessary
\begin{equation}
[D - \phi_1(x) +c ] [ D- \phi_1(x) -c]\, x=0~, \label{eq11}
\end{equation}
from where the following compatible first order ODE is obtained
\begin{equation}\label{eq12}
\dot{x}- (\phi_1(x) +c)\,x=\zeta(t)x~,
\end{equation}
or by using Eq.~(\ref{eq9c}), we obtain for this case the following Bernoulli differential equation
\begin{equation}\label{eq13} 
\dot{x}+ \tilde{c} \tan(\tilde{c} t+\delta)\, x=\phi_1(x) x~.
\end{equation}

From the expressions of $\zeta(t)$, one can see that only the one containing the tangent function is trigonometric. Therefore, taking into account that the general solutions of the Bernoulli equations can be written as logarithmic derivatives, then one may obtain isochronous (periodic) solutions of Eq.~(\ref{eq10}) only in the case of the trigonometric $\zeta(t)$. In addition, since (\ref{eq13}) is a Bernoulli equation, the factoring function $\phi_1(x)$, and so $\phi_2(x)$ too, can be of the form of a monomial plus a constant, see also the last paragraph of Section 4.

\bigskip

\section*{3. Applications}

In this section, we show how the method works in practice. We use two Li\'enard equations of the cubic oscillator type, one of them being obtained by reduction from the Sharma-Tasso-Olver equation, and Wilson's cubic-quintic oscillator case.

\ms

\noindent {\bf 3.1. Cubic anharmonic oscillator}

Let us consider the Li\'enard equation with $F(x)={\rm k}^2 x^3 + \omega^2 x$ and $G(x)=3{\rm k} x$,
\begin{equation}\label{eq22}
\ddot{x} + 3{\rm k} x \dot{x} + {\rm k}^2 x^3 + \omega^2 x= 0~,
\end{equation} 
whose unusual dynamical properties have been already noted two decades ago \cite{csl2005}, and in a radial coordinate can be related to the Lane-Emden equation in astrophysics \cite{dt1990,bluman2008}. The coefficients of this equation fulfill the relationship
\begin{equation}\label{eq22a}
F(x)=\frac{G(x)}{\sigma^2}\left[\int^x G(x')dx' +\kappa \right],
\end{equation}
for $\sigma^2=9/2$ and $\kappa=3\omega^2/2{\rm k}$. Notice that (\ref{eq22a}) is the quadrature of
Chiellini's integrability condition for Abel equations of the first kind in the dependent variable $1/\dot{x}$ to which the Li\'enard equations can be reduced, see \cite{mr2013,hlm2014}. 
Furthermore, linearizability to damped harmonic oscillators via nonlocal transformations of Li\'enard equations satisfying (\ref{eq22a}) with $\sigma \neq 0$ and
$\kappa$ arbitrary parameters was demonstrated in \cite{ks2016,s2020}. 
In addition, by the scaling $x=y/3{\rm k}$, equation (\ref{eq22}) takes the form 
\begin{equation}\label{eq22c}
\ddot{y} +  y \dot{y} + y^3/9 + \omega^2 y= 0~, 
\end{equation}
which is Liouvillian integrable since it has the form of equation (11) in \cite{sym2019} for $\alpha=2\beta^2/9$ therein.

\ms 

On the side of the factorization method, equation (\ref{eq22}) admits the simple commutative factorization
\begin{equation}\label{eq23}
\left(D + {\rm k} x +i \omega \right) \left(D + {\rm k} x - i \omega \right) x = 0~,
\end{equation}
where $\phi_1 = -{\rm k}x + i \omega$, $\phi_2 = - {\rm k} x -i \omega$, and $c=i \omega$ is an imaginary constant. Thus for this case, (\ref{eq13}) takes the form
\begin{equation} \label{eq24}
\dot{x}+\omega \tan( \omega t+\delta)x = - {\rm k} x^2~,
\end{equation}
with solution
\begin{equation}\label{eq25}
x(t) = \frac{\cos(\omega t+ \delta)}{{\rm A} + \frac{{\rm k}}{\omega} \sin(\omega t+ \delta)}~,
\end{equation}
where ${\rm A}$ and $\delta$ are arbitrary integration constant and phase, respectively. It has the following symmetries:
\begin{equation}\label{eq25a}
x(t; {\rm A},{\rm k})=-x(-t;{\rm A},{\rm k})~, \quad x(t; {\rm A},{\rm k})=x(-t;{\rm A},-{\rm k})~.
\end{equation}
If ${\rm A} \notin [-{\rm k}/\omega, +{\rm k}/\omega]$, these solutions 
are bound (regular) isochronous of period $T=2\pi/\omega$, while if ${\rm A}$ is in that range the solutions are singular isochronous.

\begin{figure}[ht!]
  \centering
 \subfigure[$\,$ Bounded and singular isochronous solutions] {\includegraphics[height=5.5cm] {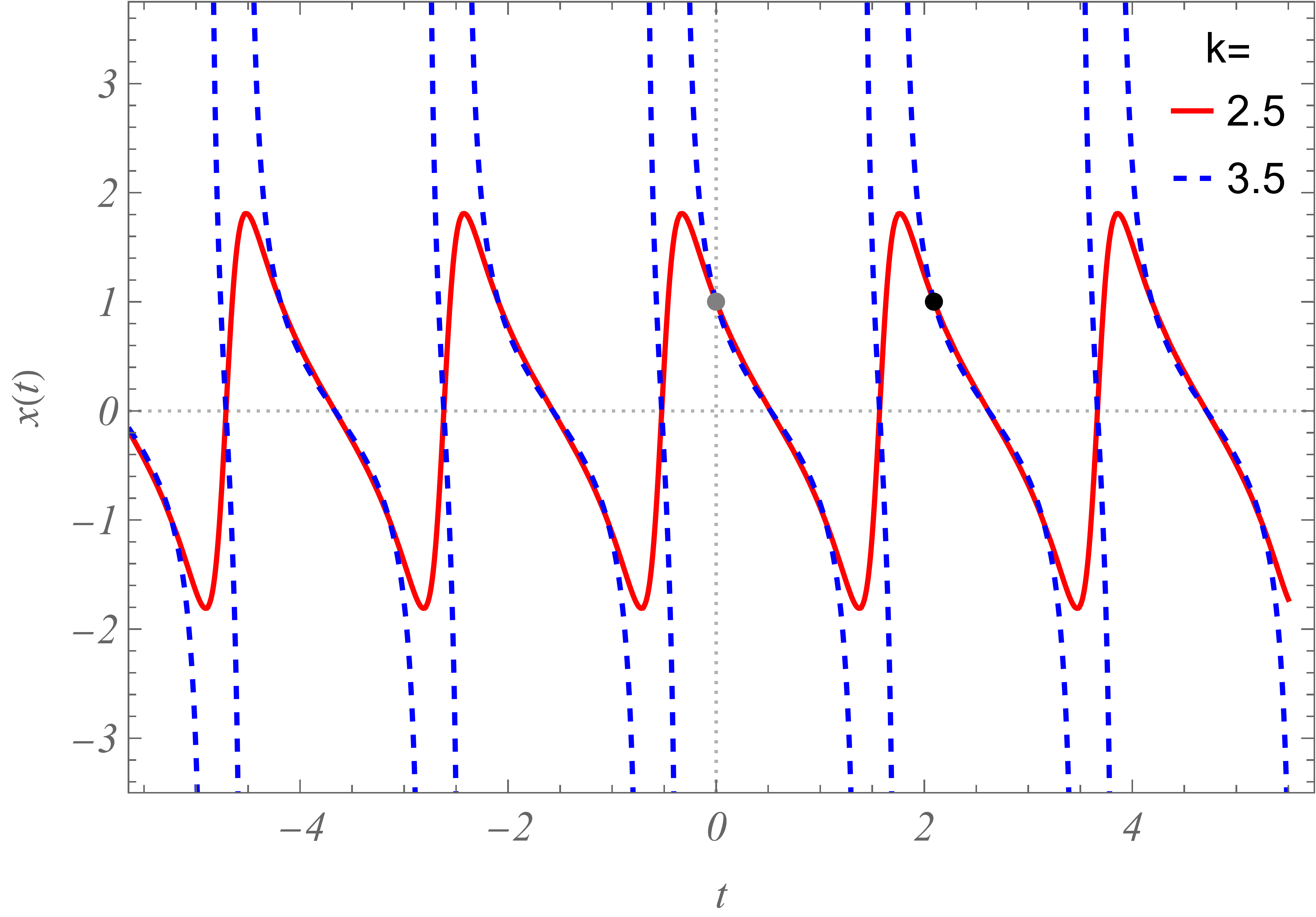}} 
    \subfigure[$\,$ The phase portraits]{ \includegraphics[height=5.5cm]  {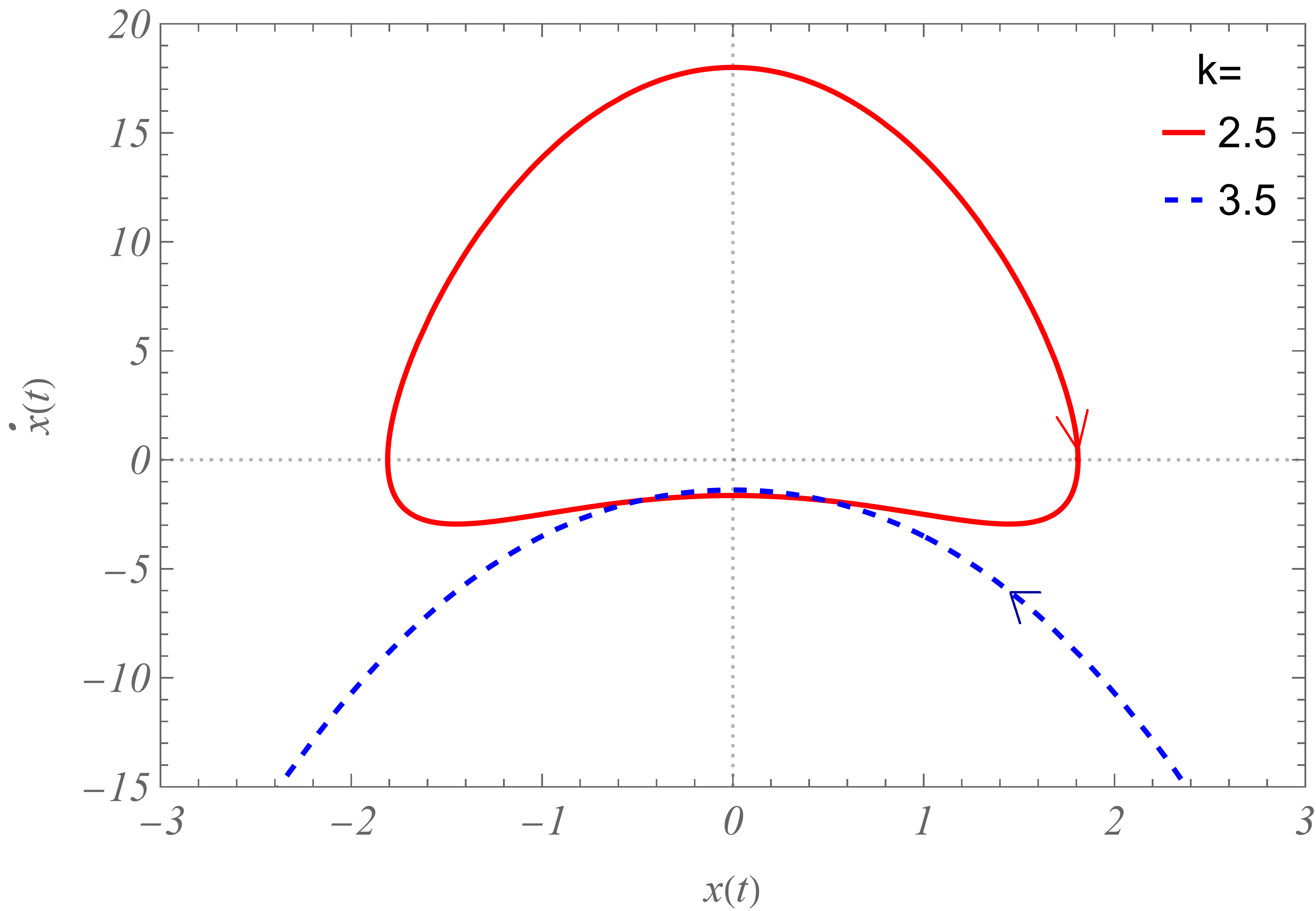}}  
   \caption{(a) Two isochronous solutions $x(t)$, of equal period as determined by the two bullets,
   of the anharmonic cubic oscillator 
   for the following values of the parameters: $|{\rm A}|=1$, $\delta=0$, $\omega=3$, and ${\rm k}=2.5$ (a bounded case) and 3.5 (a singular case).
   (b) The corresponding phase portraits.}
  \label{fig2a}
\end{figure}

\begin{figure}[ht!]
  \centering
 \subfigure[$\,$ Bounded and singular isochronous solutions] {\includegraphics[height=5.5cm] {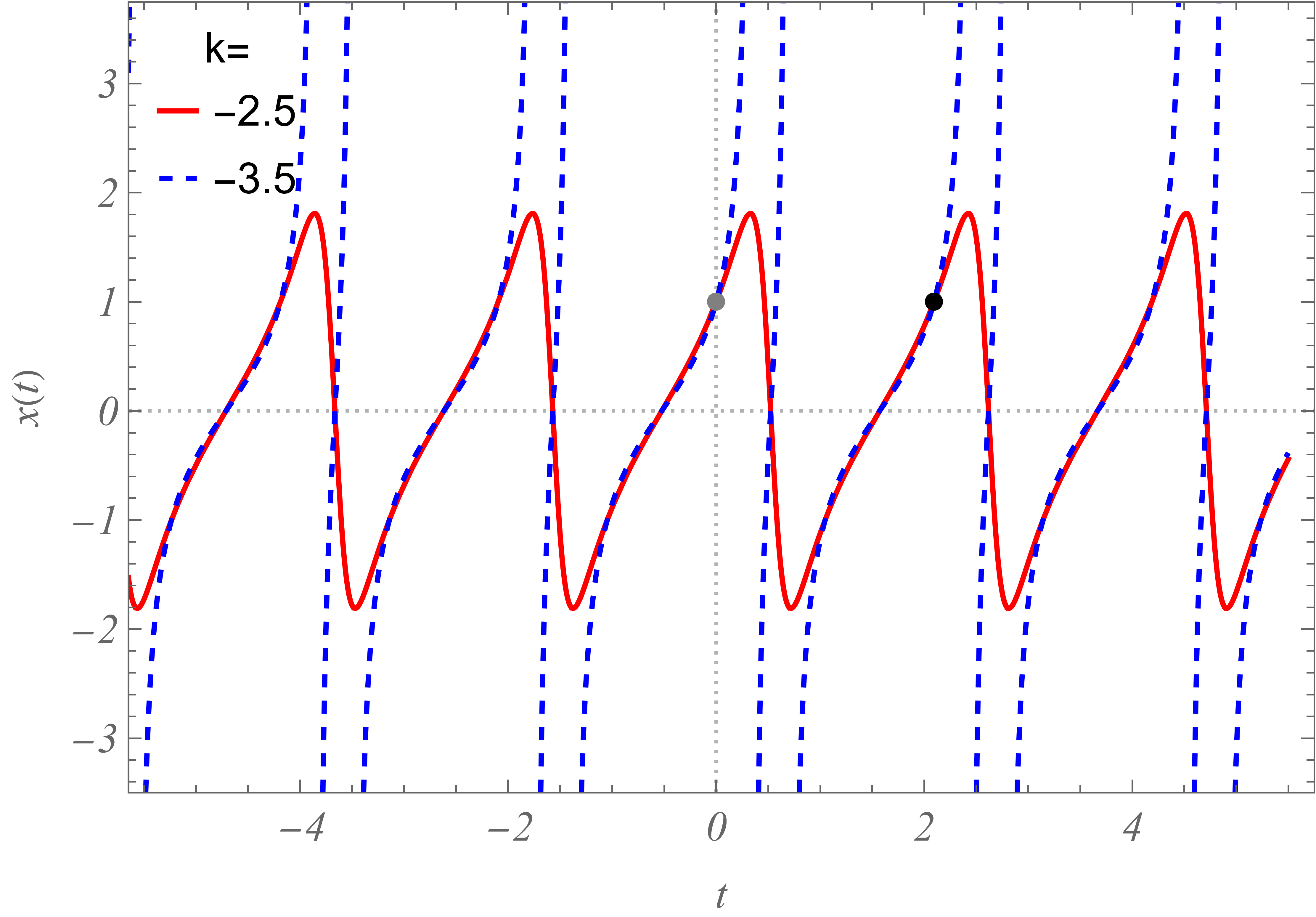}} 
    \subfigure[$\,$ The phase portraits]{ \includegraphics[height=5.5cm]  {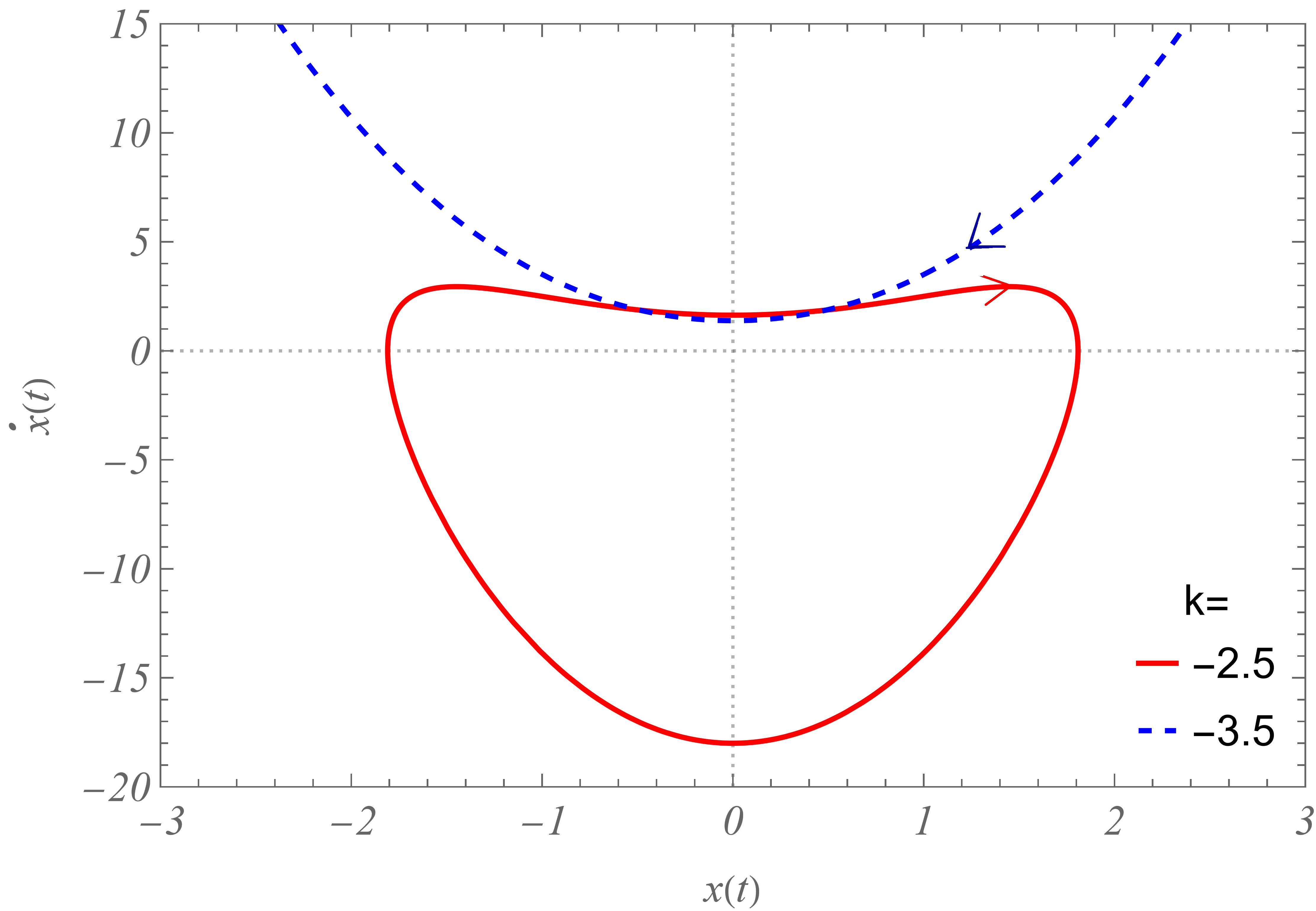}}  
   \caption{(a) Two isochronous solutions $x(t)$, of equal period as determined by the two bullets,
   of the anharmonic cubic oscillator 
   for the following values of the parameters: $|{\rm A}|=1$, $\delta=0$, $\omega=3$, and ${\rm k}=-2.5$ (a bounded case) and -3.5 (a singular case).
   (b) The corresponding phase portraits.}
  \label{fig2b}
\end{figure}

In Fig.~\ref{fig2a}(a), a regular solution and a singular one for ${\rm k}=2.5$ and ${\rm k}=3.5$, respectively, are shown together with their phase portraits
in panel (b), where one can see the open parabolic curve in the singular case that does not wrap onto the closed orbit. The same is shown in Fig.~\ref{fig2b}, 
but for the corresponding negative ${\rm k}$ parameters. In the case of the regular solutions, the closed $\dot{x}$-symmetric curves in the phase portraits 
suggest that the origin is a center as it is quite common for the Li\'enard equations with odd $F(x)$ and $G(x)$, continuous functions in some neighbourhood of the origin, and $G(x)$ of one sign for $x>0$, all these conditions according to a theorem on the existence of a centre in \cite{JordanSmith1987}. 

\ms

Regarding the singular solutions, we notice that for ${\rm A}=0$, corresponding to the centre of the interval $[-{\rm k}/\omega, +{\rm k}/\omega]$,
the solution turns into a pure periodic cotangent function of period $\pi/\omega$. For any other ${\rm A}$'s in the interval, the solutions have a slightly deformed cotangent shape with the singularities also slightly displaced to the left or to the right with respect to the singularities of the cotangent function, but are of period  $2\pi/\omega$.

\ms 

As for the arbitrary phase $\delta$, we take $\delta=0$ in all our graphs, but if one takes $\delta=\pi/2$, then in (\ref{eq24}) the cotangent function replaces the tangent function, and in (\ref{eq25}) the cosine and sine functions are interchanged.

\newpage

\ms

\noindent {\bf 3.2. Travelling wave solutions of the Sharma-Tasso-Olver equation}

The Sharma–Tasso–Olver (STO) equation is a third order nonlinear partial differential equation of the Burgers hierarchy that plays a very important role in sciences and engineering fields, for instance plasma, nonlinear optics, and quantum fields with many solution methods available in the literature \cite{ll05,ww07,jl2025}. As a member of the Burgers hierarchy, it has been linearized by generalized Cole-Hopf transformations \cite{ks2009}. 
The  STO equation is given by 
\begin{equation}\label{e1}
  u_{t}+3\alpha u^2_{x}+3\alpha u^2u_x+3\alpha uu_{xx}+\alpha u_{xxx}=0~,
\end{equation}
where the subindices define the corresponding partial derivatives. The integrability of a slightly different form of this equation obtained after scaling transformations from a Burgers equation with high order corrections has been thoroughly discussed by Kudryashov and Sinelshchikov \cite{ks2014}. Here, we show that the class of trigonometric periodic Li\'enard type solutions can be also obtained by the factorization method. The point is that the STO equation can be reduced to the following Li\'enard form,
\begin{equation}\label{e4}
  \frac{d^2U}{dz^2}+3(\epsilon+U)\frac{dU}{dz}+U^3+3\epsilon U^{2}+\left(3\epsilon^2+\frac{v}{\alpha}\right)U=0~,
\end{equation}
by firstly  
passing it to the moving reference frame 
$z=x+vt$ such that $u(x,t)=u(z)$ and integrating once to obtain the equation
\begin{equation}\label{e2}
   \frac{d^2u}{dz^2}+3u\frac{du}{dz}+u^3+\frac{v}{\alpha}u+\frac{I}{\alpha}=0~,
\end{equation}
where $I$ is an integration constant, and, secondly, by performing the shift transformation $u(z)=U(z)+\epsilon$, where the constant shift parameter $\epsilon$ satisfies the following cubic equation
\begin{equation}\label{e3}
  \epsilon^3+\frac{v}{\alpha}\epsilon+\frac{I}{\alpha}=0~.
\end{equation}
The roots of this equation are classified according to the discriminant
$$
\Delta_3=\frac{I^2}{\alpha^2}+\frac{4}{27}\left(\frac{v}{\alpha}\right)^3~.
$$
If $\Delta_3>0$, there is only one real root, and if $\Delta_3\leq 0$, all the three roots are real.

The Li\'enard equation (\ref{e4}) can be factored as follows  
\begin{equation}\label{e5}
  \left(\frac{d}{dz}+U+\frac{3}{2}\epsilon +i\sqrt{3\epsilon^2/4+v/\alpha}\right)\left(\frac{d}{dz}+U+\frac{3}{2}\epsilon -i\sqrt{3\epsilon^2/4+v/\alpha}\right)U=0~,
\end{equation}
only for the case of the real roots of the cubic equation (\ref{e3}).
For the special case when we choose the integration constant $I\equiv I_b=-2\alpha\left((1-b) \sqrt{v/\alpha}\right)^3$, where $b=0,1,2$, $\alpha >0$ and $v>0$, we have the real roots for the cubic equation given by $\epsilon_b=(1-b)\sqrt{v/\alpha}$. Therefore,
from the factorization given in Eq.~(\ref{e5}), we have that $\phi_2 -\phi_1=- i\sqrt{3\epsilon^2\alpha^2+4v\alpha}/\alpha$, which implies $\tilde{c}_b=\sqrt{v\alpha (3(1-b)^2+4)}/2\alpha$. Consequently, the general solution of Eq.~(\ref{e5}) is obtained by solving the following Bernoulli differential equation
\begin{equation}\label{e6} 
\frac{dU}{dz} + \left(\tilde{c}_b \tan( \tilde{c}_b z+\delta)+ \epsilon_b \right)U=-U^2 ~,
\end{equation}
with solution 
\begin{equation}\label{e-hw}
U_b(z) = \frac{\cos(\tilde{c}_b z +\delta)}
{Ae^{\epsilon_b\tilde{c}_b z}+
\frac{1}{\tilde{c}_b\sqrt{1+\epsilon_b^2}}\sin(\tilde{c}_b z+\delta+\varphi_b)}~, \qquad \varphi_b=\arctan(-\epsilon_b)
\end{equation}
where $\delta$ and ${\rm A}$ are integration constants. Case by case, we obtain

\noindent (i) for $b=1$, $\epsilon_1=0$ (Liouvillian integrable case),
\begin{equation}\label{e7a}
U_1(z) = \frac{\cos\left(\sqrt{\frac{v}{\alpha}} z+ \delta\right)}{{\rm A} +\sqrt{\frac{\alpha}{v}} \sin\left(\sqrt{\frac{v}{\alpha}}z+ \delta\right)}~,
\end{equation}
essentially the same as (\ref{eq25});

\noindent (ii) for $b=0$, $\epsilon_0=\sqrt{\frac{v}{\alpha}}$, 
\begin{equation}\label{e7b}
U_0(z) = \frac{\cos\left(\sqrt{\frac{7v}{4\alpha}} z+ \delta\right)}{{\rm A}e^{\sqrt{\frac{7}{4}}\frac{v}{\alpha}z} + \frac{1}{\sqrt{\frac{7}{4}(\frac{v}{\alpha}+1)}}\sin\left(\sqrt{\frac{7v}{4\alpha}}z+ \delta-\arctan\left(\sqrt{\frac{v}{\alpha}}\right)\right)}~,
\end{equation}

\begin{figure}[ht!]
  \centering
  \subfigure[$\,$ Regular and singular isochronous STO solutions]{\includegraphics[height=5.5cm] {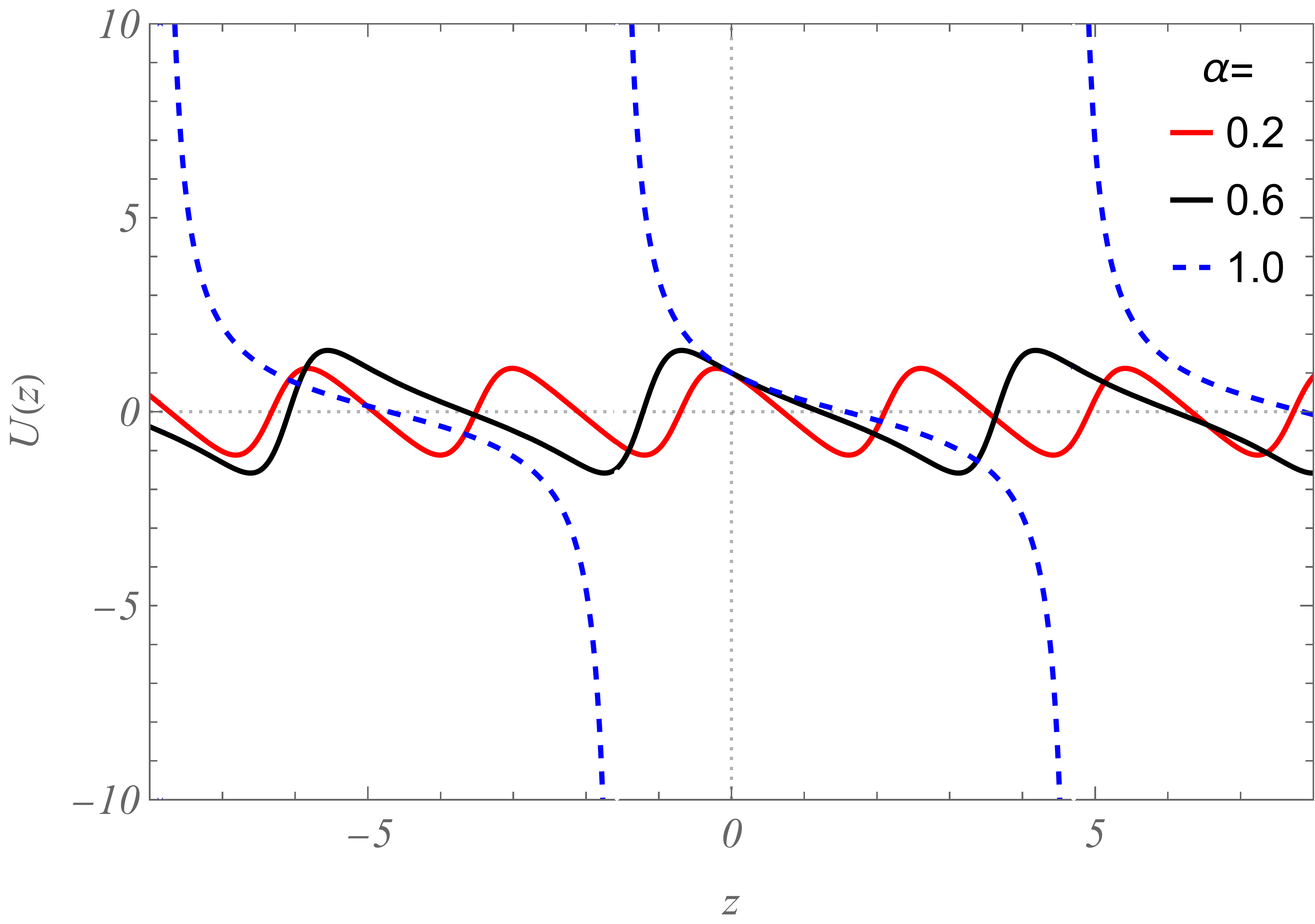}}   
  \subfigure[$\,$ Phase portraits of the solutions in (a)]{\includegraphics[height=5.5cm] {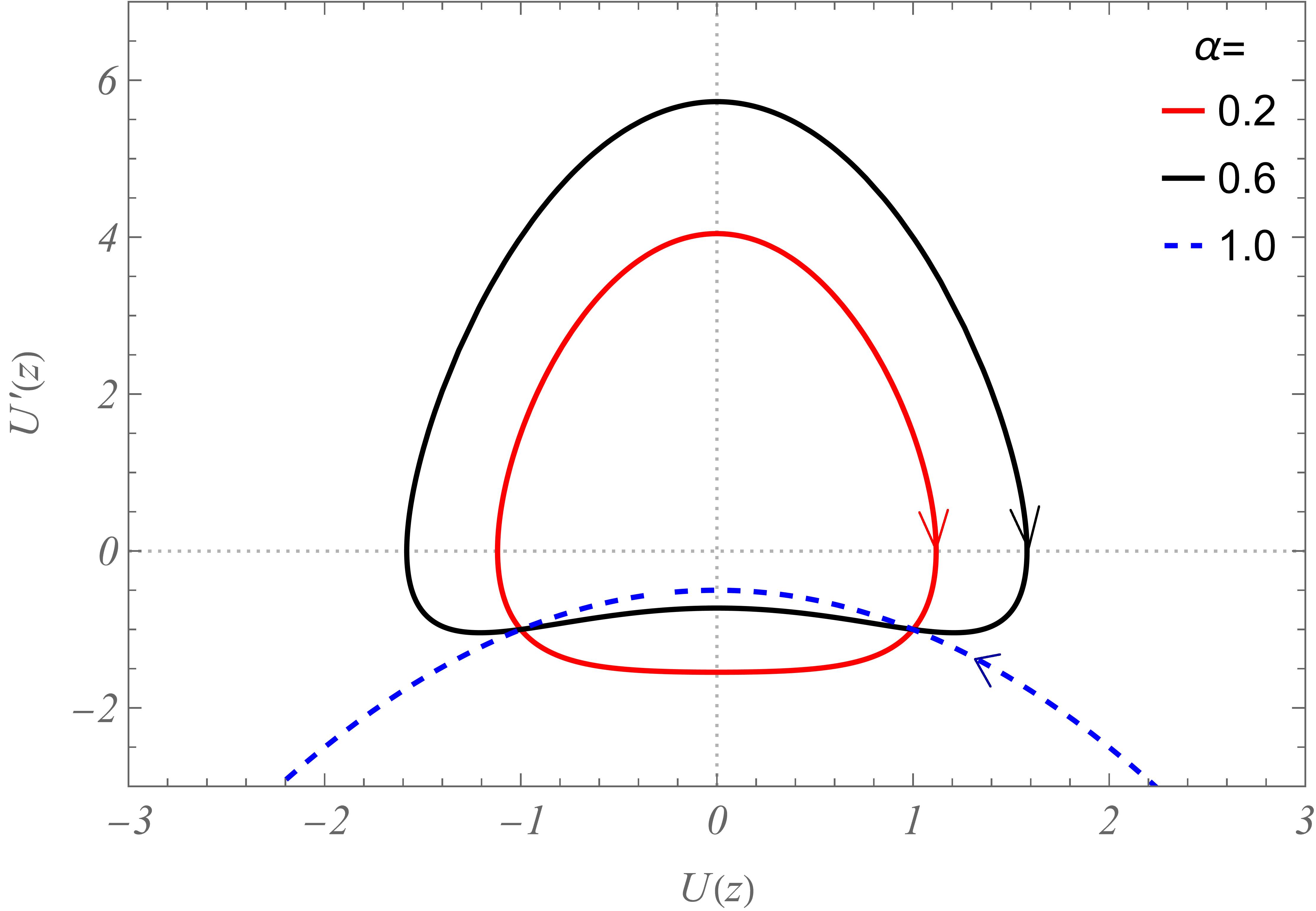}}      
   \caption{(a) Bounded and singular isochronous solutions $U(z)$ for $b=1$ ($\epsilon=0$) when they are identical to the travelling STO solutions $u(z)$.
    The other parameters are                                                          
    ${\rm A}=1$, $\delta=0$, $v=1$, and $\alpha=0.2, 0.6$, and 1. (b) The corresponding phase portraits.}
  \label{fig3}
\end{figure}

\noindent (iii) for $b=2$, $\epsilon_2=-\sqrt{\frac{v}{\alpha}}$,
\begin{equation}\label{e7c}
U_2(z) = \frac{\cos\left(\sqrt{\frac{7v}{4\alpha}} z+ \delta\right)}{{\rm A}e^{-\sqrt{\frac{7}{4}}\frac{v}{\alpha}z} +
\frac{1}{\sqrt{\frac{7}{4}(\frac{v}{\alpha}+1)}} \sin\left(\sqrt{\frac{7v}{4\alpha}}z+ \delta+\arctan\left(\sqrt{\frac{v}{\alpha}}\right)\right)}~, 
\end{equation}
 Because of the exponential functions in the denominator, full $z$-line isochronous solutions can be obtained only in the case $b=1$ for which regular isochronous solutions similar to the isochronous solution (\ref{eq25}) are obtained if ${\rm A}\notin [-\sqrt{\alpha/v},+\sqrt{\alpha/v}]$. Two such solutions are displayed
in Fig.~\ref{fig3} together with their closed orbits in phase space 
for ${\rm A}=1$, $\delta=0$, $v=1$, and two $\alpha$'s $<1$. The closed curves in the phase portraits are for the $\epsilon=0$ Li\'enard equation (\ref{e4}) which satisfies the conditions of the existence theorem of a centre at the origin \cite{JordanSmith1987}. In the same figure, one can also see the singular periodic solution for $\alpha=1$ corresponding to the right edge of the interval of singular solutions together with its parabolic orbit in phase space. The period of all these solutions is $\Lambda_{STO}=2\pi\sqrt{\alpha/v}$. 

\ms

On the other hand, if we take ${\rm A}=0$ (the centre of the singularity interval), the solutions (\ref{e7a}-\ref{e7c}) are of cotangent type,
of period $\Lambda_1=\pi\sqrt{\alpha/v}$ in the case (\ref{e7a}) and cotangent like of period $\Lambda_{0,2}=\pi\sqrt{4\alpha/7v}$ in the cases 
(\ref{e7b}) and (\ref{e7c}).

\newpage

\noindent {\bf 3.3. Wilson's Li\'enard equation}

The Wilson polynomial Li\'enard equation \cite{WilsonL} is the following particular case of the cubic-quintic Duffing-van der Pol equation
\begin{equation}\label{wpl1}
\ddot{x}+\mu(x^2-1)\dot{x}+\frac{\mu^2}{16}x^3(x^2-4)+x=0~.
\end{equation}
It was the first example of Li\'enard equation with an algebraic hyperelliptic limit cycle for $|\mu|<2$, given by the curve \cite{odani1995, qian2021}
\begin{equation}\label{wpl1a}
y^2+\frac{\mu}{2}x(x^2-4)y+(x^2-4)\bigg[\frac{\mu^2}{16}x^2(x^2-4)+1\bigg]=0~,
\end{equation}
and also admits the commutative factorization
\begin{equation}\label{wpl2}
\bigg[D+\frac{\mu}{2}\left(\frac{x^2}{2}-1\right)\pm i\frac{\sqrt{4-\mu^2}}{2}\bigg]
\bigg[D+\frac{\mu}{2}\left(\frac{x^2}{2}-1\right)\mp i\frac{\sqrt{4-\mu^2}}{2}\bigg]x=0~.
\end{equation}
For this case the constant $\tilde{c}$ is
\begin{equation}\label{wpl3}
 \tilde{c}\equiv \tilde{c}_\mu =\pm\frac{\sqrt{4-\mu^2}}{2}~, \qquad 0<|\mu|<2~
 \end{equation}
 and the corresponding Bernoulli equation (\ref{eq10}) has the trigonometric form
\begin{equation}\label{wpl4}
\dot{x}+\left(\tilde{c}_\mu\tan(\tilde{c}_\mu t+\delta)-\frac{\mu}{2}\right)x=-\frac{\mu}{4}x^3~.
\end{equation}
The general Wilson solution is then
\begin{equation}\label{wpl5}
x(t)=\pm\frac{\cos(\tilde{c}_\mu t+\delta)}
{\Big\{{\rm A}e^{-\mu t}+\frac{1}{4}+\left(\frac{\mu}{4}\right)^2\cos[2(\tilde{c}_\mu t+\delta)]+\frac{\mu\tilde{c}_\mu}{8}\sin[2(\tilde{c}_\mu t+\delta)]\Big\}^{1/2}}~,
\end{equation}
where ${\rm A}$ and $\delta$ are integration constants. To the best of our knowledge, this solution has not been mentioned in the literature yet.

\begin{figure}[ht!]
  \centering
\subfigure[$\,$]  {\includegraphics[height=5.5cm] {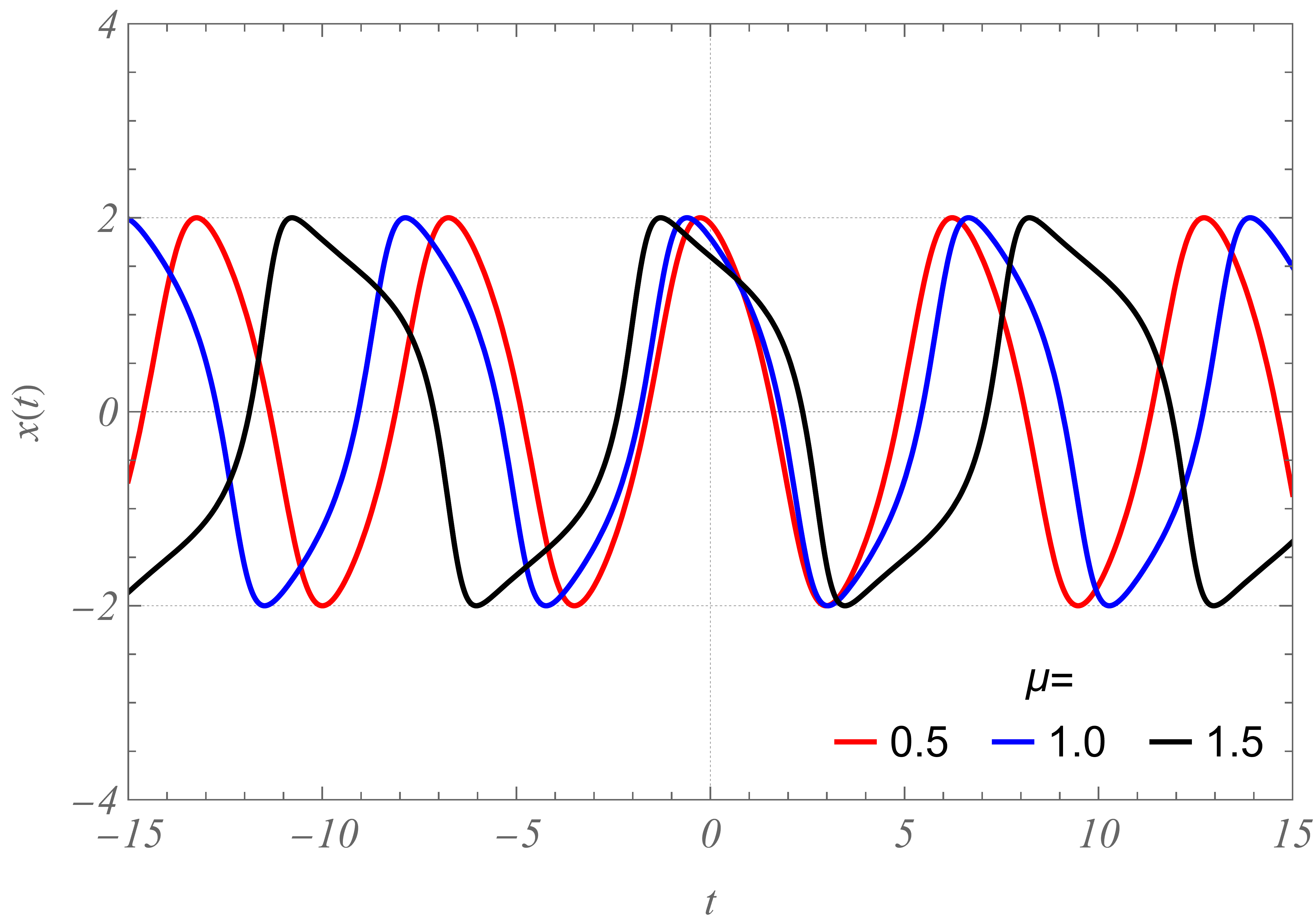}} 
\subfigure[$\,$]  {\includegraphics[height=5.5cm]  {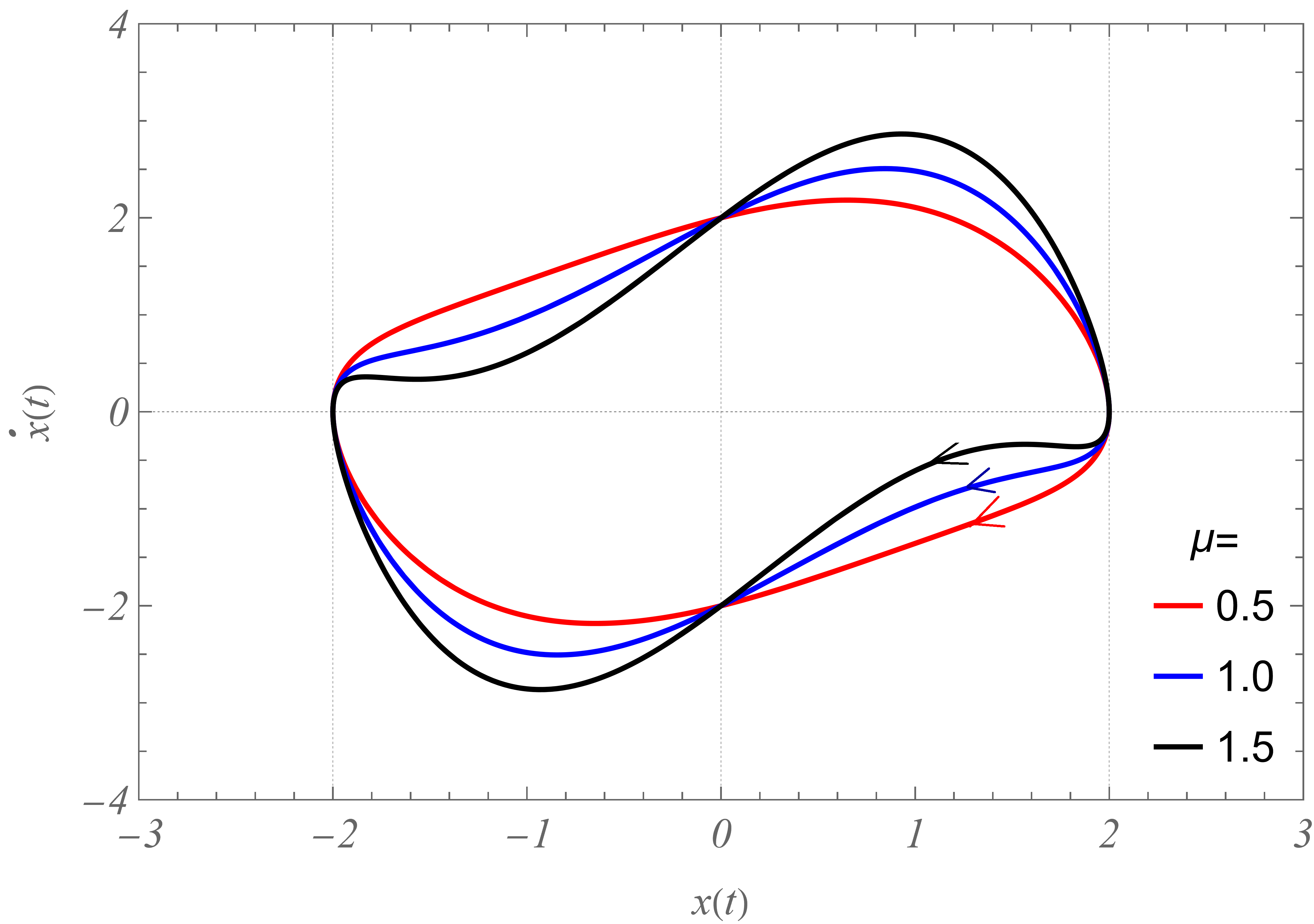}} 
   \caption{(a) Isochronous waveform solutions 
   of the Wilson polynomial Li\'enard equation obtained from (\ref{wpl5}) for ${\rm A}=0$, $\delta=0$ and three values of the parameter $\mu$; (b) their phase portraits.}
  \label{fig4}
\end{figure}

\begin{figure}[ht!]
  \centering
\subfigure[$\,$]  {\includegraphics[height=5.7cm] {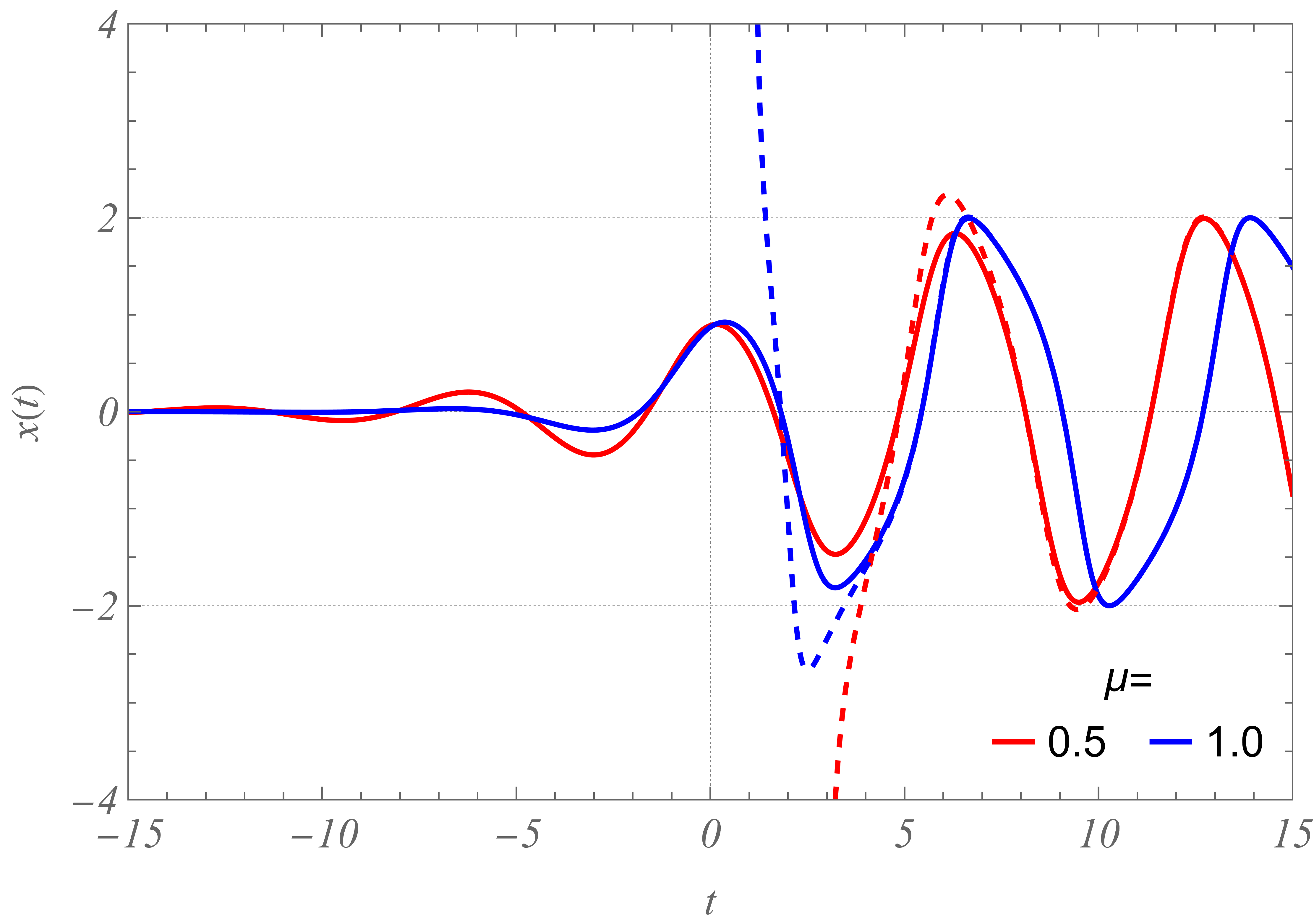}} 
\subfigure[$\,$]  { \includegraphics[height=5.7cm] {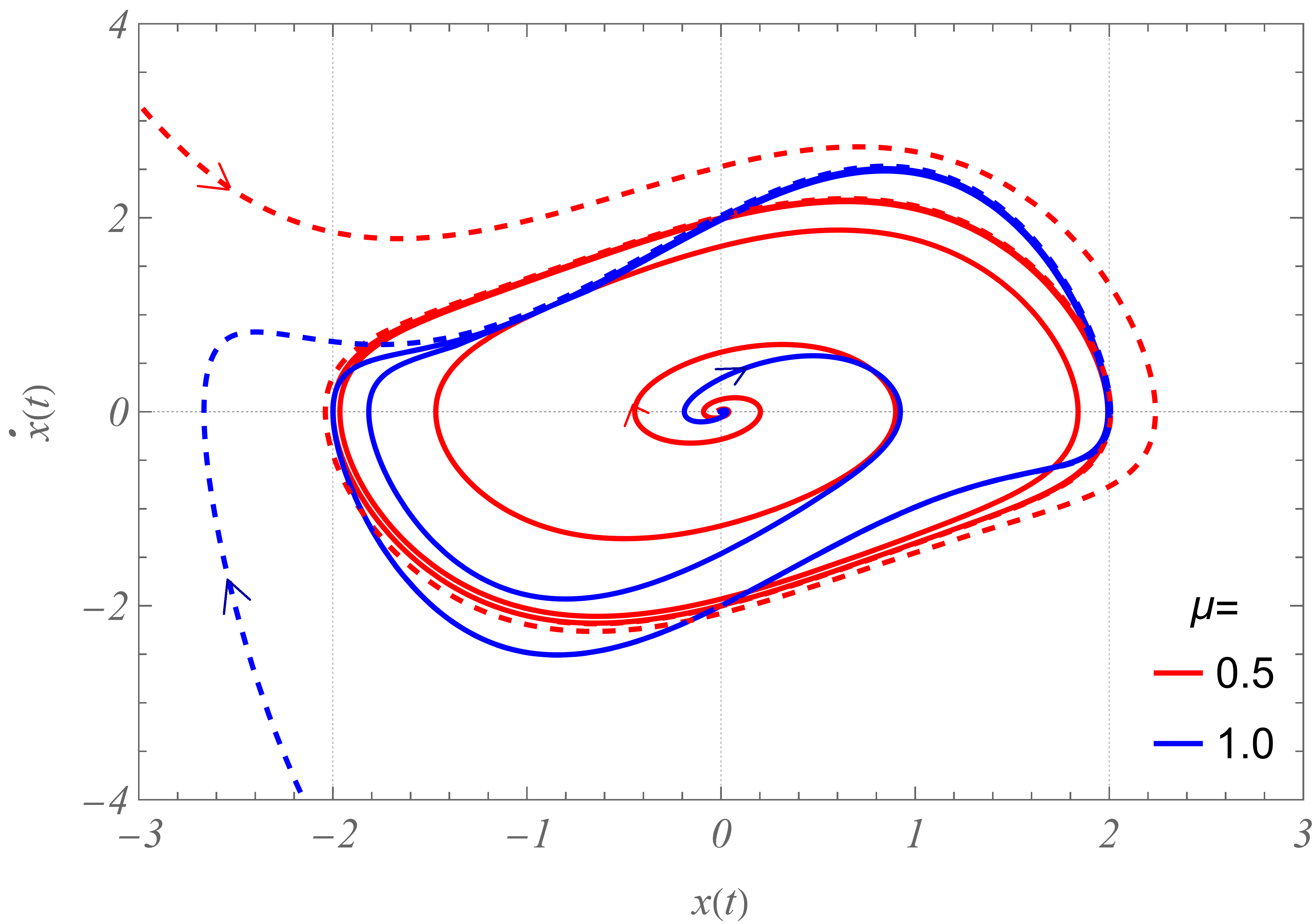}} 
   \caption{(a) Non-isochronous regular and singular Wilson waveforms from (\ref{wpl5}), unbroken lines (${\rm A}=1$) and dotted lines (${\rm A}=-1$), 
   respectively, $\delta=0$, and $\mu=$ 0.5 and 1. (b) their phase portraits.}
  \label{fig5}
\end{figure}
In Fig.~\ref{fig4}, 
one can observe isochronous Wilson solutions of type (\ref{wpl5}) and their phase portraits for ${\rm A}=0$, $\delta=0$, and three indicated values of $\mu$.
Their period is given by $T_W=2\pi/\sqrt{1-(\mu/2)^2}$.
On the other hand, in Fig.~\ref{fig5}, one can observe both non isochronous regular ${\rm A}=1$ solutions of type (\ref{wpl5}) and their phase portrait trajectories going asymptotically to the corresponding closed orbits from the inside and non isochronous singular ${\rm A}=-1$ solutions whose phase portrait trajectories approach asymptotically their closed orbit from the outside. These closed curves are limit cycles as given by (\ref{wpl1a}) for the corresponding $\mu$ parameters.

\bigskip

\section*{4. General form of polynomial Li\'enard equations that can be commutatively factored} 

In this section, we provide a theorem on the type of Li\'enard equations that are solvable through the 
commutative factorization algorithm presented in Section 2 and determine when they can also be Sundman linearizable. 

{\bf Theorem}. General solutions of Li\'enard equations of the type
\begin{equation}\label{u1}
\ddot{x}+\big[(q+2)Ax^{q}+B\big]\dot{x}+Ax^{q+1}(Ax^{q}+B)+Cx=0~, \quad q\neq 0 ~,
\end{equation}
where $A\neq 0$, $B$, and $C$ are arbitrary real parameters, are solutions of the Bernoulli differential equation
\begin{equation}\label{u4}
\dot{x}+ \bigg[\frac{B \pm \sqrt{\Delta}}{2}-\zeta_i(t)\bigg]x=-Ax^{q+1}~, \quad i= 1,2,3~, 
\end{equation}
where $ \Delta =B^2-4C$ and
\begin{eqnarray}
\zeta_1(t)&=&\frac{1}{t-t_0} \quad \textrm{for} \,\, \Delta=0~, \label{equ9a}\\
\zeta_2(t) &=& -\sqrt{\Delta}/2 \left[1- \tanh(\sqrt{\Delta}\,t/2+\delta)\right] \quad \textrm{for} \,\, \sqrt{\Delta}/2 \in \mathbb{R} \setminus \{0\}~, \label{equ9b}\\
\zeta_3(t) &=& -\sqrt{|\Delta|}/2  \left[i+\tan(\sqrt{|\Delta|}\,t/2+\delta)\right] \quad \textrm{for} \,\, \sqrt{\Delta}/2 \in \mathbb{I} \setminus \{i0\}~, \label{equ9c}
\end{eqnarray}
where $t_0$ and $\delta$ are arbitrary parameters.

\ms

The reason (proof) is that (\ref{u1}) can be factored with
\begin{equation}\label{u2}
\phi_1(x)=-\left(Ax^{q}+\frac{B}{2}\pm\frac{\sqrt{\Delta}}{2}\right)~, \quad \phi_2(x)=-\left(Ax^{q}+\frac{B}{2}\mp\frac{\sqrt{\Delta}}{2}\right)~, 
\end{equation}
which satisfy the commutative factorization condition of the form
\begin{equation}\label{u3}
\phi_2-\phi_1= \pm \sqrt{\Delta}\equiv {\rm const}.~
\end{equation}
implying the three possible cases of $\zeta(t)$ as in (\ref{eq9a})-(\ref{eq9c}), but with $c$ replaced by $\mp \sqrt{\Delta}/2$.
If $\Delta<0$, one deals with the trigonometric $\zeta_t(t)$ case leading to isochronous (periodic) solutions.

Consequently, the Li\'enard equations (\ref{u1}) have the general solutions
\begin{equation}\label{u5}
x(t)=\sqrt[q]{\frac{e^{-q\int^t Z_i(\tau)d\tau}}{{\rm K}+qA\int^t e^{-q\int^{\tau} Z_i(t')dt'}d\tau}}~,
\end{equation}
where 
${\rm K}$ is an arbitrary integration constant and $Z_i(t)=\frac{B \pm \sqrt{\Delta}}{2}-\zeta_i(t)$.
This solution can be also written in the logarithmic derivative form
\begin{equation}\label{u6}
x(t)= [qA]^{-\frac{1}{q}}\bigg[{\rm dlog}\left(\frac{{\rm K}}{qA}+\int^t e^{-q\int^{\tau} Z_i(t')dt'}\right)\bigg]^{\frac{1}{q}}~.
\end{equation}
The choice of a positive sign in front of the parenthesis for the factoring functions $\phi_1$ and $\phi_2$ in (\ref{u2}), yields a similar factorization scheme whose equivalence is obtained through $A\rightarrow -A$ and $B\rightarrow -B$, however the same form of general solution as given in Eq.~(\ref{u5}) is obtained.

If we take $q=n-1$, $n\geq 2$, a natural number, the factoring functions in (\ref{u2}) cover all the afore discussed particular cases and some others as well.
Indeed, 
$A={\rm k}, B=0, C=\omega^2$, and $n=2$ corresponds to the cubic oscillator case, $A=1, B=3\epsilon, C=3 \epsilon^2+\frac{v}{\alpha}$, and $n=2$
provides the STO case, $A=1, B=0, C=1$, and $n=2m+2$ corresponds to equation~(4) in \cite{iacono} which includes equation (\ref{eq22}) as a particular case.
The Wilson case corresponds to  $A=\mu/4, B=-\mu, C=1$, and $n=3$. Another $n=3$ case with all $A,B,C$ nonzero is a cubic-quintic Duffing-van der Pol unforced oscillator discussed in \cite{siewe}, and the list can be continued.

\ms


Interestingly, the coefficient of the first derivative of the Bernoulli integrable Li\'enard equations is similar to the monomial cases $G(x)=ay^q+b$, where $a,q\neq 0$ and $b$ are arbitrary parameters, mentioned in \cite{ks2016} as simple cases
for the application of Sundman transformations and reduction to the damped linear oscillator equation. 
The slight difference is that in (\ref{u1}) $a\equiv a_q=(q+2)A$ depends on the power of the monomial. Thus, it is clear that not all Li\'enard equations (\ref{u1}) satisfy the condition for the applicability of Sundman transformations. For equations (\ref{u1}), this condition \cite{ks2016} can be written
in the form
\begin{equation}\label{ucond}
Ax^{q+1}(Ax^q+B)+Cx=\sigma^{-2}\big[(q+2)Ax^{q}+B\big]\bigg[\frac{q+2}{q+1}Ax^{q+1}+Bx+\kappa\bigg]~,
 \end{equation}
 where $\sigma \neq 0$ and $\kappa$ are arbitrary parameters.  

It is easy to check that cases with $B=0$ allow only for $q=1$ be turned into the linear damped oscillator of unit natural frequency with the damping coefficient $\sigma=\sqrt{9/2}$. For these cases $\kappa=3C/2A$.
In our examples, the cubic oscillator in the example 3.1 and the STO case with $\epsilon=0$ are $q=1$, $B=0$ cases, and therefore they are also linearizable by Sundman transformations.

A second possibility is $\kappa=0$ that allows only for $\sigma^2=(q+2)^2/(q+1)$, $B=1$, and $C=1/\sigma^2=(q+1)/(q+2)^2$, which implies a real $0<\Delta<1$
for positive powers $q$, and a real $\Delta>1$ for negative powers ($q\neq -2$). Consequently, for these cases, one should use $\zeta_2$ from (\ref{u4}).  As examples, for $q=1$, the Li\'enard equation which is both Bernoulli and Sudman linearizable is
\begin{equation}\label{e-56}
\ddot{x}+(3Ax+1)\dot{x}+A^2x^3+Ax^2+\frac{2}{9}x=0~,
\end{equation}
with general solution
\begin{equation}\label{e-56s}
x(t)=\frac{ e^{-t/2}\cosh(t/6+\delta)}{c_1-\frac{3A}{4} e^{-t/2} \big[\sinh\left(t/6+\delta\right)+3\cosh\left(t/6+\delta\right)\big]}~,
\end{equation}
while for $q=2$, it is
\begin{equation}\label{e-57}
\ddot{x}+(4Ax^2+1)\dot{x}+A^2x^5+Ax^3+\frac{3}{16}x=0~,
\end{equation}
with general solution
\begin{equation}\label{e-57s}
x(t)=\frac{i  e^{-\frac{t-\delta}{2}}\cosh(t/4+\delta)}{\sqrt{- c_1+
\frac{8A}{3}e^{3(t/4+\delta)}\cosh^3(t/4+\delta)}}~.
\end{equation}
Solution (\ref{e-56s}) is regular for $A$ and $c_1$ of opposite signs and it is of kink type. On the other hand, 
solution (\ref{e-57s}) is real of kink type only for negative $A$ and positive $c_1$. 
%

\ms

To close this section, we mention that in a somewhat different formulation, we can say that equations of the type
\begin{equation}\label{u7}
\ddot{x}-\big[xP_x(x)+2P(x)\big]\dot{x}+\big[P^2(x)+c^2\big]x=0~,
\end{equation}
where $P(x)$ is a homogeneous polynomial of degree $n$ can be factored as
\begin{equation}\label{u8}
(D-P(x)+ic)(D-P(x)-ic)x=0~.
\end{equation}
For the examples described in this work, we have $P_k(x)=-{\rm k}x$, $P_{_{STO}}=U+3\epsilon/2$, and $P_W=\frac{\mu}{2}\left(\frac{x^2}{2}-1\right)$, respectively. These are monomial cases, 
and so they are Bernoulli integrable cases. However, when $P(x)$ is not monomial, the equation (\ref{eq13}) is integrable only in sparce cases.

\section*{5. Conclusion}

An alternative way of solving the nonlinear commutative factorization setting previously studied in \cite{getal} has been provided here that can be used to obtain periodic solutions of Li\'enard type equations. An intermediate function in the method is assumed to depend explicitly on the independent variable of the factored equation, aside to the dependent variable. This leads to a quasi-linear differential equation which can be solved and allows one to obtain the general solution of the problem for many commutative factorization settings. When a separable form in the two variables of the intermediate function is assumed, the dependent variable of the nonlinear equation satisfies a Bernoulli differential equation with parameters and nonlinear part that depend on the factorization functions. This Bernoulli differential equation is easily solved and leads to the general solution of the Li\'enard type equation. We have shown how this method works in practice in the case of three interesting mathematical physics examples of Li\'enard equations, the convective cubic anharmonic oscillator, a travelling-variable-reduced form of the Sharma-Tasso-Olver evolution equation, and the cubic-quintic Wilson's Li\'enard equation, with further insights into their solutions, including their phase portraits.
This factorization method can be used in an efficient manner when polynomial factoring functions that differ by a constant can be found. The general form of the Li\'enard equation for an arbitrary polynomial degree $n$ that can be factored in this way is provided as well as the conditions to be also approached by non-local Sundman transformations.

\bigskip
\bigskip

\noindent {\bf CRediT authorship contribution statement}\\
\noindent {\bf G.~Gonz\'alez}: Conceptualization, Writing – original draft.

\noindent {\bf O. Cornejo-P\'erez}: Methodology, Validation.

\noindent {\bf J. de la Cruz}: Methodology, Formal analysis.

\noindent {\bf H.C.~Rosu}: Formal analysis, Writing - review and editing.

\bigskip
\bigskip

\noindent {\bf Declaration of competing interests}\\
The authors declare that they have no known competing financial interests or personal relationships
that could have appeared to influence the work reported in this paper.

\bigskip
\bigskip

\noindent {\bf Data availability statement}\\
No data have been used in this paper.

\bigskip
\bigskip

\noindent {\bf Acknowledgments}\\
The authors wish to thank the referees for their remarks that led to substantial improvements of this work.
G.G. would like to acknowledge support by the program C\'atedras Conahcyt through project 1757 and from project A1-S-43579 of SEP-CONAHCYT Ciencia B\'asica and Laboratorio Nacional de Ciencia y Tecnolog\'{\i}a de Terahertz.

\end{document}